\documentclass{acmart}

\usepackage{mathrsfs}

\newcommand{\I}{\mathrm{I}}

\newcommand{\trace}{\mathrm{tr}}
\newcommand{\rank}{\mathrm{rank}}
\newcommand{\range}{\mathrm{range}}
\newcommand{\Stiefel}{\mathrm{St}}
\newcommand{\Grassman}{\mathrm{Gr}}
\newcommand{\Orth}{\mathrm{Orth}}
\newcommand{\GL}{\mathrm{GL}}

\newcommand{\Bigo}{\mathscr O}

\newcommand{\RN}{\ensuremath{\mathbb{R}}}

\newcommand{\detFlow}{{\tt det-flow}}
\newcommand{\trFlow}{{\tt trace-flow}}

\begin{document}

\title{PCA by Determinant Optimization has no Spurious Local Optima}

\author{Raphael A.\ Hauser}
\affiliation{
  \institution{University of Oxford and Alan Turing Institute}
  }
\author{Armin Eftekhari}
\affiliation{%
 \institution{Alan Turing Institute and University of Edinburgh}
}

\author{Heinrich F.\ Matzinger}
\affiliation{
  \institution{Georgia Institute of Technology}
}

%



\begin{abstract}
{Principal component analysis (PCA)  is an indispensable tool in many learning tasks that finds the best linear representation for data. Classically, principal components of a dataset are interpreted as the directions that preserve most of its ``energy'', an interpretation that is theoretically underpinned by the celebrated Eckart-Young-Mirsky Theorem. There are yet other ways of interpreting PCA that are rarely exploited in practice, largely because it is not known how to reliably solve the corresponding non-convex optimisation programs.  In this paper, we consider one such interpretation of principal components as the directions that preserve most of the ``volume'' of the dataset. Our main contribution is a theorem that shows that the corresponding non-convex program has no spurious local optima. We apply a number of solvers for empirical confirmation.}

\end{abstract}
\maketitle

\section{Introduction}\label{introduction}


Let $A\in\mathbb{R}^{m\times n}$ be a data matrix whose rows correspond to $m$ different items, and columns to $n$ different features. 
For $p\le n$ and $X\in\mathbb{R}^{n\times p}$ with orthonormal columns, rows of $AX\in\mathbb{R}^{m\times p}$ correspond to the original data vectors projected onto $\range(X)$,  the column span of $X$. In particular, the new data matrix $AX$ has reduced dimension $p$, while the number $m$ of items is unchanged. 
Successful dimensionality reduction is at the heart of classification, regression and other learning tasks that often suffer from the ``curse of dimensionality'', where having a small number of training samples in relation to the data dimension (namely, $m\ll n$) typically leads to overfitting  \cite{hastie2013elements}. 
Principal component analysis (PCA) is one of the oldest dimensionality 
reduction techniques going back to the work of Pearson \cite{pearson} and Hotelling \cite{hotelling}, which were motivated by the observation that often data is approximately located in a lower-dimensional subspace. PCA identifies this subspace by finding a suitable matrix $X$ that retains in $AX$ as much as possible of the energy of $A$. The optimal $X$ is called the \emph{loading matrix}. 
The columns of the loading matrix also yield important structural information about the data by identifying groups of variables that occur with jointly positive or jointly negative weights. In various application domains such groups of variables indicate patterns of functional dependencies, for example, genes that are jointly upregulated \cite{alter4}.
%
%
More formally, assume throughout this paper that  the data matrix $A$ is mean-centred, namely $\sum_{i=1}^m a_i = 0$ where $a_i\in\mathbb{R}^n$ is the $i$-th row of $A$. Consider $p\le n$ and let $\mathbb{R}^{n\times p}_p$ be the space of full column-rank $n\times p$ matrices and consider the {\em trace inflation function}
\begin{eqnarray}
g_{t}:\RN^{n\times p}_p&\rightarrow&\RN\nonumber\\
X&\mapsto&
\frac{\| AX\|_F^2}{\|X\|_F^2} = 
\frac{\trace(X^*A^*AX)}{\trace(X^*X)}\label{gtrace}
\end{eqnarray}
and the program
\begin{equation}\label{pca by trace}
\arg\max\left\{g_{t}(X):\;X\in\Stiefel(n,p)\right\}.
\end{equation}
Above, $\|\cdot\|_F$ and $\mbox{tr}(\cdot )$ return the Frobenius norm and  trace of a matrix, respectively, and $A^*$ is the transpose of matrix $A$. 
With $p\le n$,  $\Stiefel(n,p)$ above denotes the the Stiefel manifold, namely the set of all $n\times p$ matrices with orthonormal columns.  
Program \eqref{pca by trace} performs PCA on the data matrix $A$, and it is a consequence of the Eckart-Young-Mirsky Theorem \cite{eckart, mirsky} that a Stiefel matrix $X\in\Stiefel(n,p)$ is a global optimizer of Program \eqref{pca by trace} if and only if it is a $p$-leading right-singular factor $V_p\in \mathbb{R}^{n\times p}$ of $A$, namely if and only if $V_p$ consists of  the right singular vectors of $A$ corresponding to its $p$ largest singular values.  
We remark that when $X\in \Stiefel(n,p)$, the denominator in the definition of $g_t(X)$ in \eqref{gtrace} is constant and serves the purely cosmetic role of highlighting the similarities with an alternative program to compute principal components of $A$ that we will propose below, a program that is unconstrained and does not require $X$ to have orthogonal columns. 
The interpretation of PCA as a tool for dimensionality reduction suggests that it should suffice to merely find a matrix 
$X\in\RN^{n\times p}_p$ whose columns span the desired optimal subspace obtained by optimization over the Grassmannian $\Grassman(n,p)$, defined as the set of 
$p$-dimensional subspaces of $\RN^n$. However, program \eqref{pca by trace} does not fit that bill, for it is easy to see that $g_t(X)$ is not invariant under a change of basis for $\range(X)$. 

\section{Main Result and Prior Art}

In analogy to the function $g_t$ in \eqref{gtrace}, we define a {\em volume inflation function} by 
\begin{eqnarray}
g_{d}:\RN^{n\times p}_p&\rightarrow&\RN\nonumber\\
X&\mapsto&\frac{\det(X^*A^*AX)}{\det(X^*X)},\label{gdet}
\end{eqnarray}
where $\det$ stands for determinant. 
Unlike  $g_t$, note that $g_d$ is  invariant under a general change of basis. Indeed, consider 
$X\in\RN^{n\times p}_p$, $\Theta\in\GL(p)$, and $Y=X\Theta$, where $\GL(p)$ is the general linear group, i.e., the set of nonsingular $p\times p$ matrices. Then we have that  
\begin{equation}
g_{d}(Y)= \frac{\det(\Theta)^2\det(X^*A^*AX)}{\det(\Theta)^2\det(X^*X)}
=g_{d}(X).
\label{invariance of det}
\end{equation}
It follows that $g_d$ is naturally defined on the Grassmannian $\Grassman(n,p)$. 
In analogy to \eqref{pca by trace}, consider the program 
\begin{equation}\label{pca by det}
\arg\max\left\{g_{d}(X):\;X\in\RN^{n\times p}_p\right\}.
\end{equation}
Note that Program \eqref{pca by det} is unconstrained because $\RN^{n\times p}_p$ is an open subset of $\RN^{n\times p}$ with nonempty relative interior, and that, having chosen $X\in\RN^{n\times p}_p$, we explicitly exclude degenerate points.  

Program \eqref{pca by det} is a good model for dimensionality reduction because $X^*A^*AX$ is the sample covariance of the dimensionality-reduced data matrix $AX$, and $\det(X^*A^*AX)^{\frac{1}{2}}$ is the volume of the smallest box that fits around the ellipsoid $\{y\in\RN^p:\,y^*(X^*A^*AX)^{-1}y\leq 1.96^2\}$, which closely approximates the smallest ellipsoid that encloses all of the projected (dimensionality-reduced) data points $(a_i X)^*$, where $a_i$ is the feature vector of the $i$-th item, that is the $i$-th row of matrix $A$. Thus, Program \eqref{pca by det} may be thought of as maximizing  the volume of point cloud of dimensionality reduced items. Program \eqref{pca by trace} in contrast maximizes 
the length of the diagonal of the aforementioned box.  
Thus, Programs \eqref{pca by trace} and \eqref{pca by det} 
have different geometric justifications. Our main result is an extension of the Eckart-Young-Mirsky Theorem to Program \eqref{pca by det}, showing that the latter yields the same notion of PCA of $A$ as Program \eqref{pca by trace}:

\begin{theorem}\label{det optimality}
The following statements hold true: 
\begin{itemize}
\item[i) ] Any $p$-leading right-singular factor $V_p\in\mathbb{R}^{n\times p}$ of $A$ is a global maximizer of Program (\ref{pca by det}). 
\item[ii) ] For any global maximizer $\widetilde{X}$ of Program (\ref{pca by det}),  there exists a $p$-leading right-singular factor $V_p$ of $A$ such that $\range(\tilde{X})=\range(V_p)$. 
\item[iii) ] Program (\ref{pca by det}) does not have any spurious local maximisers, namely any local maximizer of Program  (\ref{pca by det}) is also a global maximizer.
\end{itemize}
\end{theorem}

Due to the characterization of $X$ in part ii),  the singular value decomposition of 
$A\widetilde{X}=\widetilde{U}\widetilde{S}\widetilde{V}^*$, which can be computed in merely 
$\Bigo(m p^2)$ time, yields $V_p=X\widetilde{V}$ as a $p$-leading right-singular factor of $A$, the diagonal coefficients 
of $\widetilde{S}$ as the $p$ leading singular values of $A$, and the first $p$ columns of $\widetilde{U}$ as the corresponding 
$p$-leading left singular factor of $A$.  
Parts i) and ii) of Theorem \ref{det optimality} are not new, as they can be  easily proven using interlacing properties of singular values 
or via the Cauchy-Binet Formula 
. What is new, however, is Part iii), which is crucial in rendering the computation of a PCA via Program (\ref{pca by det}) practical: this property guarantees that any  locally convergent descent algorithm from the standard literature (gradient descent, coordinate descent, block coordinate descent, trust-region methods, line search descent methods, cubic regularisation methods, etc.) applied to the negative objective function of (\ref{pca by det}) are automatically globally convergent to a correct PCA  irrespective of the starting point used. Empirical confirmation of this property is provided in our numerical section where we compare the convergence of gradient descent, line-search gradient descent and accelerated gradient descent algorithms applied to Programs (\ref{pca by trace}) and (\ref{pca by det}). 
To prove Property iii), we develop an exact differential geometric characterization of all KKT points of Program (\ref{pca by det}), including the case with coalescing singular values, and this also yields a novel, purely geometric proof of Parts i) and ii) of Theorem \ref{det optimality}. In contrast to earlier proofs of these parts, our analysis does not depend on the columns of $X$ being mutually orthogonal.
%
%
A string of conceptually related papers has appeared in the recent literature, all of which show the nonexistence of spurious local optima for different problems for the purposes of understanding the geometry and performance of iterative descent algorithms on nonconvex optimization problems \cite{ge-lee-ma, sun2015nonconvex, GHJY15, SQW15, 
BNS16
}. Our theory fits nicely into this growing literature, although it is based on different techniques. Among this line of research papers, the work of  \cite{li} is most closely related to ours, as they consider functions $X\mapsto\phi(X^*BX)$ with $B$ positive semidefinite and $\phi$ convex, while our objective function $X\mapsto \log(g_{d}(X))$ is a difference of two such functions. 

\section{Technical Details}\label{details}

Before proving Theorem \ref{det optimality}, a word on the notation is in order. 
We denote the singular value decomposition of a matrix $A\in\RN^{m\times n}$ of rank $k$ by $A=USV^*$, 
where $U\in\Orth(m)$, $V\in\Orth(n)$, and $S\in\RN^{m\times n}$ is a diagonal matrix with the ordered singular values $\sigma_1\geq 
\sigma_2\geq \dots\geq\sigma_k>0$  of $A$ on the diagonal and all other values  equal to zero.  If $k<n$ 
we complete the spectrum by setting $\sigma_{k+1}=\dots=\sigma_n=0$. The {\em thin SVD} of $A$ will be 
denoted by $A=U_k S_k V_k^*$, where $U_k\in\Stiefel(m,k)$ consists of the first $k$ columns of $U$, 
$S_k$ of the top left $k\times k$ block of $S$, and $V_k\in\Stiefel(n,k)$ of the first $k$ columns of $V$. 
For any index set $J\subseteq\{1,\dots,n\}$, 
we write $V_J$ for the matrix composed of the columns $\{v_j:\,j\in J\}$, and $S_J$ for the submatrix 
$[s_{ij}]$ of $S$ with indices $j\in J$, $i\in J\cap\{1,\dots,m\}$, and $U_J$ for the matrix consisting of the columns 
$\{u_i :\,i\in J\cap\{1,\dots,m\}$, leaving the original ordering intact in all cases. We also write 
$A_J=U_J S_J V_J^*$ for the low-rank 
approximation of $A$ corresponding to singular vectors 
with index in $J$. In the special case $J=\{1,\dots,p\}$ that corresponds to the $p$-leading part SVD of $A$ 
we use the shorthand notation $U_p,S_p,V_p$ and $A_p$ for $U_J,S_J,V_J$ and $A_J$ respectively. 

Let us now turn to the proof of Theorem \ref{det optimality}. 
Assuming $\rank(A)\geq p$, so that there exist $X$ for which $\det(X^*A^*AX)\neq 0$, we may instead maximize 
the function $f_{d}(X):=\ln g_{d}(X)$ 
and reformulate Problem \eqref{pca by det} in the following equivalent form, 
\begin{equation}\label{pca by det prime}
X^*=\arg\max\left\{f_{d}(X):\;  X\in\RN^{n\times p}_p\right\}
\end{equation}
which is unconstrained, because $\RN^{n\times p}_p$ is an open full-dimensional domain in $\RN^{n\times p}$. The stationary points of this model are therefore characterized by $\nabla f_{d}(X)=0$. 
Using a Taylor expansion, it is easy to see that the gradient of $f_{d}$ is given by 
\begin{equation}\label{lemma1}
\nabla f_{d}(X)=2X(X^*X)^{-1} -2A^*AX (X^*A^*AX)^{-1}
\end{equation}
and behaves as follows under a change of basis $Y=X\Theta$, where $\Theta\in\GL(p)$, 
\begin{equation}\label{lemma3}
\nabla f_{d}(Y)=\nabla f_{d}(X)\Theta^{-*}. 
\end{equation}

\begin{lemma}\label{lemma4}
Let $A\in\RN^{m\times n}$ be of rank $k\geq p$, $X\in\RN^{n\times p}_p$ and $AX\in\RN^{m\times p}_p$. Then  
$\nabla f_d(X)=0$ if and only if there exists a set of indices $J=\{i_1,\dots,i_{\ell}\}\subseteq\{1,\dots,k\}$ such that 
$|\{\sigma_{i_1},\dots,\sigma_{i_{\ell}}\}|\leq p$ (not counting multiplicities), and  $\range(X)\subseteq\range(V_J)$. 
\end{lemma}

\paragraph{Proof:} 
By virtue of Equation \eqref{lemma3}, 
$\nabla f_{d}(X)=0$ if and only if $\nabla f_{d}(X\Theta)=0$ for any $\Theta\in\GL(p)$, 
and since it is also the case that $\range(X)\subseteq\range(V_J)$ if and only if $\range(X\Theta)\subseteq\range(V_J)$, 
we may in fact assume without loss of generality that $X\in\Stiefel(n,p)$, so that $X^*X=\I_p$. 
Let $A=U_k S_k V_k^*$ be 
the thin SVD of $A$, which is an exact factorization because $\rank(A)=k$, so that 
$\sigma_1\geq\dots\geq\sigma_k>0=\sigma_{k+1}=\dots=\sigma_{n}$. By \eqref{lemma1}, we then have 
\begin{eqnarray}
\nabla f_{d}(X)&=&0
\;\Leftrightarrow\; X = A^*AX(X^*A^*AX)^{-1}
\;\Leftrightarrow\; (XX^*) V_k S_k^2 V_k^*X=V_k S_k^2 V_k^*X
\;\Leftrightarrow\;\range(V_k S_k^2 V_k^*X)\subseteq\range(X)\nonumber\\
&\Leftrightarrow& \range(V_k S_k^2 V_k^* X)=\range(X),\label{rank condition}\\
&\Leftrightarrow& \exists\,Z\in\GL(p)\;\mathrm{ s.t. } V_k S_k^2 V_k^* X=XZ
\;\Leftrightarrow\; S_k^2 V_k^*X = V_k^*XZ
\;\Leftrightarrow\; (X^*V_k) S_k^2 = Z^*(X^*V_k),\nonumber
\end{eqnarray}
where Equivalence \eqref{rank condition} follows from $\rank(V_k S_k^2 V_k^*X)=p$. Denoting the column vectors of $X^*V_k$ by $(w_1,\dots,w_k)$, we have $w_j\in\RN^p$ with 
\begin{equation}\label{eigenw}
\sigma_j^2 w_j = Z^*w_j,\quad (j=1,\dots,k). 
\end{equation}
Since $Z^*$ has at most $p$ different eigenvalues, it must be the case that $|\{\sigma_j:\,w_j\neq 0\}|\leq p$ (not counting multiplicities), and then $AX\in\RN^{m\times p}_p$ implies that we have $X=V_k V_k^* X=V_J V_J^* X$, where $J=\{j:\,w_j\neq 0\}$ is the index set whose existence is claimed in the lemma.\hfill $\square$

\begin{lemma}\label{lemma5}
Let $A\in\RN^{m\times n}$ be of rank $k\geq p$, $X\in\RN^{n\times p}_p$ such that $AX\in\RN^{m\times p}_p$ and 
$\nabla f_{d}(X)=0$, and let $[w_1,\dots,w_k]=X^*V$. 
If $\{j_1,\dots,j_s\}\subseteq\{1,\dots,k\}$ is an index set such that for any strict subset $J\subset\{j_1,\dots,j_s\}$ the vectors $\{w_j:\,j\in J\}$ are linearly independent but the set $\{w_{j_i}:\,i=1,\dots,s\}$ is linearly dependent, then $\sigma_{j_1}=\dots=\sigma_{j_s}$. 
\end{lemma}

\paragraph{Proof:}
By the assumptions of the theorem, there exist unique scalars $\lambda_2,\dots,\lambda_s$, all nonzero, such that 
$w_{j_1}=\sum_{i=2}^s \lambda_i w_{j_i}$, and multiplying this equation with the matrix $Z^*$ constructed in the proof 
of Lemma \ref{lemma4} yields $\sigma_{j_1}^2 w_{j_1}=\sum_{i=2}^s \sigma_{j_i}^2\lambda_i w_{j_i}$. By the uniqueness and non-nullity of the scalars $\lambda_i$, and by strict positivity of $\sigma_{j_1}$, comparison of the two equations yields $\sigma_{j_i}=\sigma_{j_1}$ $(i=2,\dots,s)$. 
\hfill$\square$

Lemma \ref{lemma5} shows that proper linear dependence among the vectors $w_j$ can only occur within subsets of vectors that correspond to the same eigenvalue of the matrix $Z^*$ from the proof of Lemma \ref{lemma4}. We also note that since $Z^*\in\GL(p)$, it has at most $p$ linearly independent eigenvectors. This motivates the following definition, 
in which $A\in\RN^{m\times n}$ is a matrix of rank $k\geq p$, $X\in\RN^{n\times p}_p$ such that 
$AX\in\RN^{m\times p}_p$. Let $[w_1,\dots,w_k]:=X^*V_k$, and $J:=\{j:\,w_j\neq 0\}$. Further, let
$\varsigma_1>\dots>\varsigma_{\rho}>0$ be an enumeration of the nonzero singular values 
$\sigma_1\geq\dots\geq\sigma_k$ of $A$ without repetition. 
For any index set ${\mathcal J}\subseteq\{1,\dots,k\}$ let us write $\mathcal{I}({\mathcal J})=\{i:\,\varsigma_i\in\{\sigma_j:\,j\in J\}\}$, and conversely,  for any index set ${\mathcal I}\subseteq\{1,\dots,\rho\}$ we write ${\mathcal J}({\mathcal I})=\{j:\,\sigma_j\in{\mathcal I}\}$. For $i=1,\dots,\rho$, we call $m(i):=|{\mathcal J}(\{i\})|$ is the multiplicity of $\varsigma_i$ in $V$, $m_{\lambda}(i,p):=\max\left(0,p-\sum_{\ell=1}^{i-1} m(\ell)\right)$ the $p$-leading multiplicity of $\varsigma_i$ in $V$, $m_{\tau}(i,p):=\max\left(0,p-\sum_{\ell=i+1}^{k} m(\ell)\right)$ the $p$-trailing multiplicity of $\varsigma_i$ in $V$, and $m(X,i):=\rank\{w_j:\,j\in J_i\}$  the multiplicity of $\varsigma_i$ in $X$, where $J_i=J\cap{\mathcal J}(\{i\})$.
For an intuitive interpretation, note that $m(i)$ is the number of right-singular vectors of $A$ associated with the singular value $\sigma_i$, $m_{\lambda}(i,p)$ the number of column vectors in any $p$-leading right-singular factor $V_p$ of $A$ associated with the singular value $\sigma_i$, and $m(X,i)=\dim(\range(X)\cap\mathrm{V}_i)$ where $\mathrm{V}_i:=\range\{v_j:\,\sigma_j=\varsigma_i\}$ is the right-singular space of $A$ associated with the 
singular value $\varsigma_i$. 

\begin{corollary}\label{characterization of stationarity}
$X$ is a stationary point of Program (\ref{pca by det prime}) if and only if 
$m(X,i)=\dim(\Pi_{\mathrm{V}_i}\range(X))$ for $i=1,\dots,\rho$,   
where $\Pi_{\mathrm{V}_i}$ denotes the orthogonal projection into the right-singular space of $A$ associated with the 
singular value $\varsigma_i$. 
\end{corollary}

\paragraph{Proof:}
Note that since $\range(X)\cap\mathrm{V}_i\subseteq\Pi_{\mathrm{V}_i}\range(X)$, we have 
\begin{equation}\label{multineq}
m(X,i)\leq\dim(\Pi_{\mathrm{V}_i}\range(X))
\end{equation}
for all $X\in\RN^{n\times p}$ and $i=1,\dots,\rho$. Thus, the corollary 
says in fact that $X$ is a stationary point of Program \eqref{pca by det prime} if and only if $\range(X)\cap\mathrm{V}_i=\Pi_{\mathrm{V}_i}\range(X)$ for all $i$. By Lemma \ref{lemma4}, 
$\nabla f_d(X)=0$ if and only if $\range(X)\subseteq\range(V_J)$. Assuming that the latter condition is 
satisfied, we have $\range(X)\cap\mathrm{V}_i=\{0\}=\Pi_{\mathrm{V}_i}\range(X)$ for all 
$i\notin{\mathcal I}(J)$, while the condition also implies step \eqref{oho} 
in the following sequence of equalities, 
\begin{eqnarray}
\sum_{i=1}^\rho m(X,i)&=& \sum_{i=1}^{\rho} \rank(X^*V_{J_{i}})
\;=\;\rank(X^*V_J)\label{from lemma5}\\
&=&\rank(X^*V_k)\label{oho}\\
&=&p,\label{from rank condition}
\end{eqnarray}
where \eqref{from lemma5} follows from Lemma \ref{lemma5} and \eqref{from rank condition} from 
$AX\in\RN^{n\times p}_p$. But this sequence of equalities can only hold if \eqref{multineq} holds as an equality 
for all $(i=1,\dots,\rho)$, and hence, $\range(X)\cap\mathrm{V}_i=\Pi_{\mathrm{V}_i}\range(X)$ for all 
$i\in{\mathcal I}(J\cap{\mathcal J}(\{1,\dots,\rho\}))$. Conversely, if 
$\range(X)\cap\mathrm{V}_i=\Pi_{\mathrm{V}_i}\range(X)$ for all $i$, then 
\begin{eqnarray}
\range(X)&=&\range(VV^*X)\;=\;\range(V_kV_k^*X)\label{rrrank}\\
&=&\range(V_J V_J^* X)\label{from J}\\
&\subseteq&\range(V_J),\nonumber
\end{eqnarray}
as claimed, where \eqref{rrrank} follows from $AX\in\RN^{n\times p}_p$, and \eqref{from J} from the 
definition of $J$, by which $X^*v_j=0$ for $j\notin J$.\hfill$\square$

\begin{theorem}\label{optimal value}
If $X$ is a stationary point of Program (\ref{pca by det prime}), then 
\begin{itemize}
\item[i) ] $\sum_{i=1}^{\rho} m(X,i) = p$,
\item[ii) ] $f_{d}(X)=\sum_{i=1}^{\rho} 2 m(X,i)\ln\varsigma_i$.
\end{itemize}
\end{theorem}

\paragraph{Proof:}
In assuming that $X$ is a stationary point of Program \eqref{pca by det prime}, we make the assumption that  $A\in\RN^{m\times n}$ is a matrix of rank $k\geq p$, and $X\in\RN^{n\times p}_p$ such that 
$AX\in\RN^{m\times p}_p$ and $\nabla f_{d}(X)=0$, so that all of the above lemmas apply. 
Part i) is then an equivalent reformulation of Corollary \ref{characterization of stationarity}, by virtue of 
Equation \eqref{from rank condition}. For Part ii), 
we may assume without loss of generality that $X^*X=\I_p$, for the same reasons as in the 
proof of Lemma \ref{lemma4}, so that 
$f_{d}(X)=\ln\det(X^*V_J S_J^2 V_J^* X)$. Then there exists $\Theta\in\Orth(p)$ such that 
$X\Theta = [Y_1,\dots,Y_{\rho}]$, where some of the blocks may have zero columns, and where each block 
$Y_i=[y^i_1,\dots,y^i_{m(X,i)}]\in\Stiefel\left(n,m(X,i)\right)$ satisfies $\range(Y_i)\subseteq\range(V_{J_i})$ so 
that $Y_i^* V_{J_j}=0$ $(j\neq i)$ and $V_{J_i}V_{J_i}^*Y_i=Y_i$. We have 
\begin{eqnarray*}
f_{d}(X) &=&\ln\det\left(\Theta^*X^*V_J S_J^2 V_J^* X\Theta\right)
\;=\;\ln\det\left(\left[\oplus_{i=1}^{\rho}Y_i^* V_{J_i}\right]\left[\oplus_{i=1}^{\rho}
\varsigma_i^2\I_{m(X,i)}\right]\cdot 
\left[\oplus_{i=1}^{\rho}V_{J_i}^* Y_i\right]\right)\\
&=&\ln\prod_{i=1}^{\rho} \det\left(Y_i^* V_{J_i}\left(\varsigma_i^2\I_{m(X,i)}\right) V_{J_i}^* Y_i\right)
\;=\;\ln\prod_{i=1}^{\rho}\varsigma_i^{2m(X,i)}\det\left[\left(Y_i^* V_{J_i} V_{J_i}^*\right)
\left(V_{J_i}^* V_{J_i} Y_i\right)\right]
\;=\;\sum_{i=1}^{\rho} 2 m(X,i)\ln\varsigma_i,
\end{eqnarray*}
where $\oplus$ denotes the composition of a block-diagonal matrix out of its constituent blocks.\hfill 
$\square$

\begin{theorem}\label{saddle points}
Let $X$ be a stationary point of Program (\ref{pca by det prime}). The following hold true:
\begin{itemize}
\item[i) ] $X$ is a global maximizer of $f_{d}$ if and only if $m(X,i)=m_{\lambda}(i,p)$ for all $(i=1,\dots,\rho)$.
\item[ii) ] When $\rank(A)=n$, then $X$ is a global minimizer of $f_{d}$ if and only if $m(X,i)=m_{\tau}(i,p)$ for all 
$(i=1,\dots,\rho)$.
\item[iii) ] In all other cases $X$ is a saddle point of $f_{d}$.
\end{itemize}
\end{theorem}

\paragraph{Proof:}
It suffices to prove Claim iii), as Claims i) and ii) will then follow immediately from Theorem \ref{optimal value}.iii) and the fact that $f_{d}$ is $C^{\infty}$ on $\{X\in\RN^{n\times p}:\,\rank(AX)=p\}$. 
Using the same transformation $Y=X\Theta$ as in the proof of Theorem \ref{optimal value}, we may replace $X$ by $Y$ and assume without loss of generality that $X$ consists of $p$ different columns of $V$, ordered by non-increasing corresponding singular values, $X=[v_{j_1},\dots,v_{j_p}]$, $\sigma_{j_1}\geq\dots\geq\sigma_{j_p}$, and since we are neither in the case of Claims i) or ii), there exist $\mu,\eta\in\{1,\dots,k\}\setminus\{j_1,\dots,j_p\}$ such that $\sigma_{\mu}>\sigma_{j_p}$ and $\sigma_{\eta}<\sigma_{j_1}$. Now let us write $c=\cos\theta$, $s=\sin\theta$, and consider the matrix 
$X(\theta)=\left[\begin{array}{cccc}v_{\eta}&v_{j_1}&\dots&v_{j_p}\end{array}\right]
\left[\begin{array}{cc}\tiny{\left[\begin{array}{cc}c&s\\-s&c\end{array}\right]}&0\\
0&\I_{p-1}\end{array}\right]\left[\begin{array}{cc}0\\ \I_p\end{array}\right]$. 
We have $X(0)=X$, $X(\theta)\in\Stiefel(n,p)$ for all $\theta$, and 
\begin{eqnarray*}
& &\hspace{-1cm}\det\left(X(\theta)^*A^*AX(\theta)\right)
\;=\;\det \left(X(\theta)^*VS^2V^*X(\theta)\right)\nonumber\\
&=&\det\tiny{\Bigg(\left[\begin{array}{cc}0&\I_p\end{array}\right]
\left[\begin{array}{cc}\tiny{\left[\begin{array}{cc}c&-s\\s&c\end{array}\right]}&0\\0&\I_{p-1}\end{array}\right] }
\tiny{\left[\begin{array}{cccc}\sigma_{\eta}^2&&&0\\
&\sigma_{j_1}^2&&\\
&&\ddots&\\
0&&&\sigma_{j_p}^2\end{array}\right]
\left[\begin{array}{cc}\tiny{\left[\begin{array}{cc}c&s\\-s&c\end{array}\right]}&0\\0&\I_{p-1}\end{array}\right]
\left[\begin{array}{c}0\\ \I_p\end{array}\right]\Bigg)}\nonumber\\
&=&(\sigma_{\eta}^2 s^2 + \sigma_{j_1}^2 c^2)\sigma_{j_2}^2\dots\sigma_{j_p}^2. 
\end{eqnarray*}
Thus, writing $\phi(\theta)=f_{d}(X(\theta))$, we find 
$\phi''(\theta)  = -\left(\frac{2\sigma^2_{\eta}\sin\theta\cos\theta-2\sigma^2_{j_1}\cos\theta\sin\theta}{\sigma_{\eta}^2\sin^2\theta+\sigma_{j_1}^2\cos^2\theta}\right)^2 
+\frac{2(\sigma^2_{\eta}-\sigma^2_{j_1})(\cos^2\theta-\sin^2\theta)}{\sigma_{\eta}^2\sin^2\theta+\sigma_{j_1}^2\cos^2\theta}$, and hence, $\phi''(0)= 2(\sigma^2_{\eta}-\sigma^2_{j_1})/\sigma^2_{j_1}<0$. Analogously, using 
$X(\theta) =\left[\begin{array}{cccc}v_{j_1}&\dots&v_{j_p}&v_{\mu}\end{array}\right]
\left[\begin{array}{cc}\I_{p-1}&0\\
0&\tiny{\left[\begin{array}{cc}c&s\\
-s&c\end{array}\right]}\end{array}\right]
\left[\begin{array}{c}\I_p\\0\end{array}\right]$
yields $\phi''(0)= 2(\sigma^2_{\mu}-\sigma^2_{j_p})/\sigma^2_{j_p}>0$. The Hessian of $f_{d}$ at $X$ has thus both directions of positive and negative curvature, and $X$ is a saddle point, as claimed.\hfill
$\square$\\

Theorem \ref{det optimality} is now an easy corollary of Theorem \ref{saddle points}:


\paragraph{Proof of Theorem \ref{det optimality}}
Part i) follows from the trivial fact that $m(V_p,i)=m_{\lambda}(i,p)$, $(i=1,\dots,\rho)$ and the ``if part'' of 
Theorem \ref{saddle points}.i), whereas Part ii) follows from Corollary \ref{characterization of stationarity} and the ``only if part'' of Theorem \ref{saddle points}.i).\hfill$\square$\\

\section{Experimental Verification of Results}\label{experiments}

\begin{figure}
\begin{center}
\scalebox{0.05}{\includegraphics{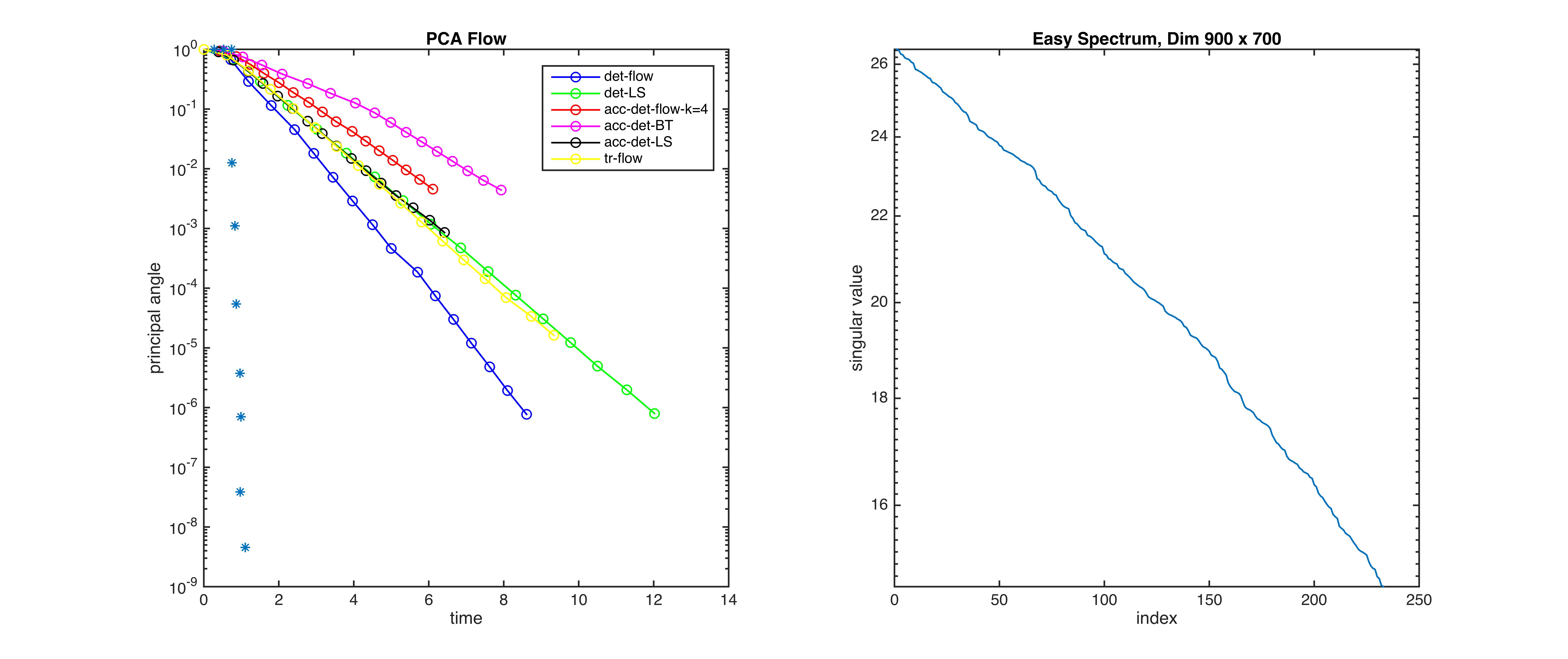}}
\caption{Computing the principal components of a $900\times 700$  random  matrix an ``easy'' spectrum using \trFlow, \detFlow,~and its variants. The spectrum is displayed on the right and the error versus time for various algorithms are plotted on the left.  \label{easy spec}}
\end{center}
\end{figure}

\begin{figure}
\begin{center}
\scalebox{0.05}{\includegraphics{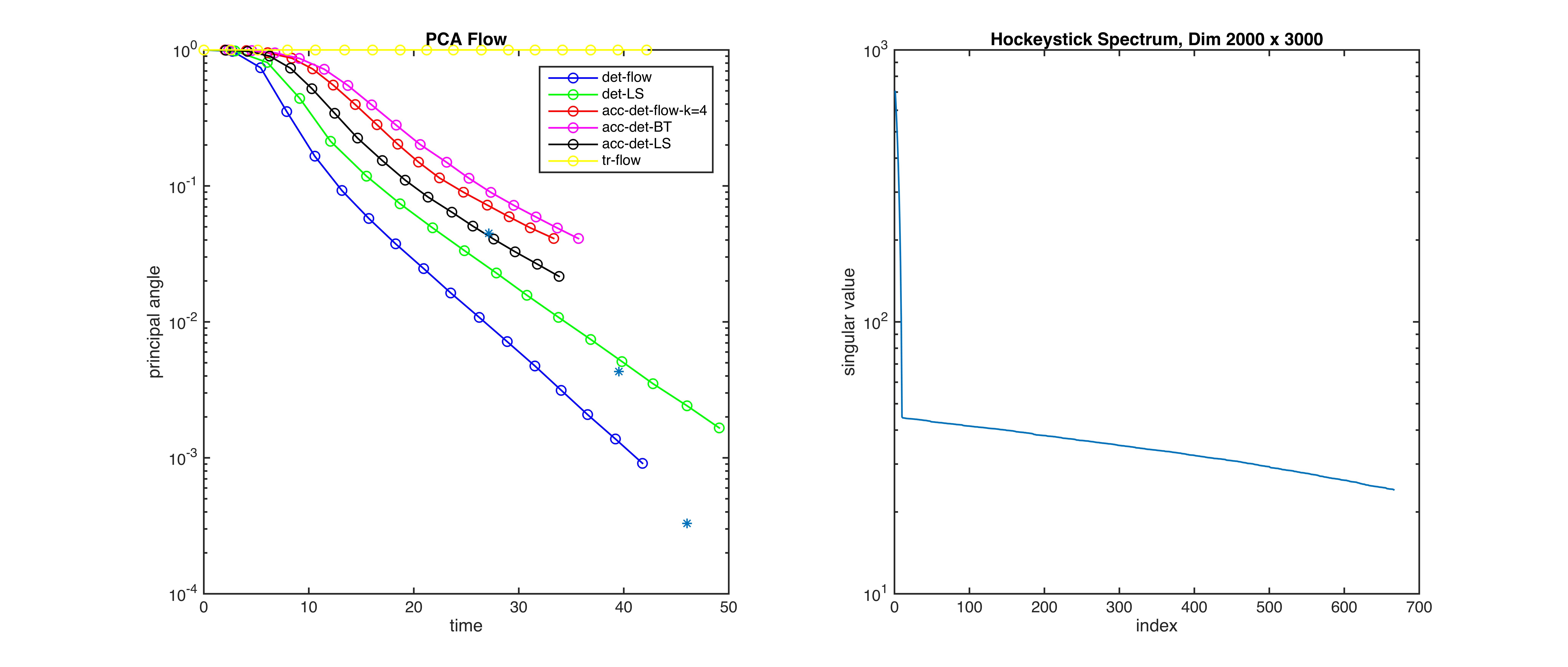}}
\caption{Similar to Figure \ref{easy spec} but for a $2000\times 3000$ random matrix with a ``hockey-stick'' spectrum. \label{hockey}}
\end{center}
\end{figure}

The theory of Section \ref{details} gives deep insight into the geometry of low-rank factorizations. Part iii) of Theorem \ref{det optimality} furthermore shows that Model (\ref{pca by det}), or equivalently, (\ref{pca by det prime}), is of practical use in conjunction with any standard descent algorithm that is known to be locally convergent, as the Theorem then implies that all such algorithms are also globally convergent irrespective of the starting point. In this section we conduct a series of experiments that give empirical confirmation of this finding. For the purposes of these experiments we chose but a few of the simplest and easiest to implement algorithms: steepest ascent with constant step size, steepest ascent with variable step size, and steepest ascent with Scieur-d'Aspremont-Bach acceleration. 
In order to get a qualitative understanding of the behavior of Model (\ref{pca by det prime}), we compare it against numerical results obtained by applying the same algorithms to model 
\begin{equation}\label{pca by trace prime}
\arg\max\left\{f_{t}(X):=\log g_t(X):\;X\in\Stiefel(n,p)\right\},
\end{equation}
which is often used in the sparse PCA context \cite{absil}. As we will see below, determinant maximization has significant advantages over trace maximization. 
It is important to keep in mind that we are comparing two optimization {\em models}, rather than {\em algorithms}, as any number of known optimization algorithms could have been applied instead of the steepest ascent variants we discuss, and a comprehensive comparative study of algorithms is beyond the scope of this paper. However, based on our extensive numerical experiments we believe that the qualitative differences between the determinant and trace models are robust under the choice of algorithm for their solution. 
Recall that the steepest ascent algorithm applied to the unconstrained maximization of a $C^1$ function $f:\RN^{n\times p}_p \rightarrow\RN$ is an iterative procedure with updating rule 
\begin{equation}\label{steepest ascent updates}
X_{i+1}= X_i +\eta \nabla f(X_i),
\end{equation}
where $\eta>0$ is a steplength multiplier. When $\eta$ is chosen constant, a good choice is $\eta=1/\Lambda$, where $\Lambda$ is a global Lipschitz constant of $\nabla f$ if such a constant is know. If $\Lambda$ is unknown, it has to be estimated numerically. The case where $\eta$ is chosen as the exact minimizer of the one-dimensional optimization problem $\min_{\eta} f(X_i+\eta\nabla f(X_i))$ is called steepest ascent with exact line-search. Since this is computationally expensive, practical line-search methods have been developed to find an approximately optimal $\eta$ via computationally inexpensive interval searches. One of the best standard choices of such a method is a line-search with {\em Wolfe Conditions}. A standard stopping criterion for unconstrained iterative optimization methods is $\|\nabla f(X)\|<\epsilon$, that is, the algorithm is stopped when the objective gradient falls below a chosen threshold. All of the above is standard textbook material, see e.g.\ \cite{nocedal2006numerical}. 
For ease of reference, we call any algorithm applied to Model (\ref{pca by det prime}) as a {\em determinant flow}, while referring to any algorithm applied to Program (\ref{pca by trace prime}) as a {\em trace flow}. Note that while (\ref{pca by det prime}) is an unconstrained optimization problem, (\ref{pca by trace prime})  is constrained by the requirement that $X$ be a Stiefel matrix of size $n\times p$. To get around this problem, \cite{absil} apply a standard steepest ascent update (\ref{steepest ascent updates}) followed by re-orthogonalization of the columns of $X$, which is obtained by applying a thin QR factorization to $X$ \cite{golub} and replacing $X$ by $Q$. This method is guaranteed to work generically, and we adopted it for all implementations of determinant and trace flows,  in order to make running times more comparable. 
In our experiments we observed that numerical estimates of the local Lipschitz constant of $\nabla f_d(X)$ was nearly constant over the entire domain of the determinant flow, hence it appears that numerically this model is globally Lipschitz constant and amenable to steepest ascent with constant step size, an algorithm we denote by \detFlow. In contrast, the trace flow numerically rotates the initial iterates $X$ very fast for most starting points, resulting in large local Lipschitz constants for $\nabla f_t(X)$ initially, and becoming small  asymptotically. The use of a constant step size $\eta$ is not appropriate for this model, as it would force a choice that renders all but the initial few step sizes too short. We therefore implemented the trace flow with a practical line search based on Wolfe Conditions \cite{nocedal2006numerical} and refer to this algorithm as \trFlow. 
To give the determinant flow the same chance at running with a variable step size, we implemented several other variants: {\tt det-LS} is the steepest ascent method with the same practical line-search, but applied to the determinant flow model rather than the trace flow. A third variant,  {\tt acc-det-flow} is an implementation of the acceleration scheme of \cite{NIPS20166267} applied to the determinant flow. This scheme computes $k+1$ consecutive updates (\ref{steepest ascent updates}) and then uses the gradients $\nabla f_d(X_k)$ to compute an improved single update direction $\Delta_i$ applied to $X_i$, 
\begin{equation*}
X_{i+1}=X_{i}+\xi \Delta_i. 
\end{equation*}
We used constant values $\eta=1/\Lambda$ for the $k+1$ updates (\ref{steepest ascent updates}), where $\Lambda$ is a numerical estimate of the Lipschitz constant of $\nabla f_d$, as in algorithm \detFlow. We tested three variants of this method: {\tt acc-det-flow-k=4} is a basic version with $\xi=1$ throughout, {\tt acc-det-LS} is a version with $k=4$ and $\xi$ chosen by a practical line search based on Wolfe Conditions, and {\tt acc-det-BT} is a version with $k=4$ and $\Delta_i$ computed via a backtracking method proposed in \cite{NIPS20166267}. The technical details of the acceleration method goes beyond the limited space available to this paper, but our implementation is fully reproducible by referring to \cite{NIPS20166267}. Furthermore, we will make our code publicly available upon acceptance of our paper. 
We compared \trFlow, \detFlow and all its variants on the following 
two random matrix models: 
Left and right singular vectors of matrices were drawn from the uniform (Haar) distribution on the corresponding spaces. For example, the left singular vectors of a random $2000\times 3000$ matrix were drawn from the uniform distribution on the orthogonal group $\Orth(2000)$, whereas its right singular vectors were drawn from the uniform distribution on the Stiefel manifold $\Stiefel(3000,2000)$. For the singular values, the first model uses an ``easy'' spectrum that has a fairly flat scree plot. 
In the second model, the singular spectrum exhibits a ``hockey-stick'' shape. 
%
We computed $p=15$ right-leading singular vectors of matrices generated from these three models and with different matrix dimensions. The error is expressed as the principal angle \cite{golub} $\arcsin(\|V_pV_p^*\widehat{X}-\widehat{X}\|)$ between $\range(\widehat{X})$ and the ground truth $\range(V_p)$. Here, $\|\cdot\|$ is the spectral norm, $\widehat{X}$ is the output of the algorithm invoked, and $V_p$ is a $p$-leading right-singular factor of the random data matrix at hand. Error of all algorithms versus time are plotted in Figures \ref{easy spec} and \ref{hockey}. 
Our code was implemented in MATLAB (Release R2015a), and all numerical experiments were carried out on a MacBook Air equipped with a 1.4 GHz Intel Core i5 processor with 8 GB of memory. While the purpose of this section is primarily to illustrate that iterative algorithms applied to the determinant flow model (\ref{pca by det prime}) have the regular convergence behavior predicted by the theory of Section \ref{details} and to give a qualitative comparison with Program (\ref{pca by trace prime}), it is also interesting for comparison to run the compiled Fortran code of the LAPACK implementation of Lanczos' Method on the same input data. The Lanczos Method (see e.g.\ \cite{golub}) was the algorithm of choice for the computation of $p$-leading part SVDs for many decades, although it has now been superseeded. An asterisk is plotted if all $p$ right-leading singular vectors were computed to the accuracy specified on the vertical axis in the time shown on the  horizontal axis. In some of the experiments the error is shown as $10^0$, and this means that not all $p$ 
vectors were computed to the required accuracy. This may not be reflective of the maximum error among the $p$ vectors, as the Lanczos code only outputs the vectors it could compute to the required accuracy. In contrast, the error of the outputs of our own implementation are reflective of the maximum error, namely the principal angle $\arcsin(\|V_pV_p^*\widehat{X}-\widehat{X}\|)$ described earlier. 
%
On the difficult spectra the \detFlow~algorithm significantly and consistently outperformed the \trFlow, the latter being competitive only for random matrices with an easy spectrum. 
We remark that although {\tt det-LS, acc-det-flow-k=4, acc-det-BT} and {\tt acc-det-LS} all have lower iteration complexity than 
 {\tt det-flow}, that is, these algorithms converge to a given target accuracy in fewer iterations, the cost per iteration is much higher due to the evaluation of a variable step size and/or the accelerated update direction $\Delta_i$, so that the overall clock time of these algorithms is slower than {\tt det-flow}. In the case of the trace flow however, the line search is very cheap, because evaluating a trace takes only $O(p)$ time. Algorithm {\tt tr-flow} has by far the lowest cost per iteration, but its iteration complexity is much higher than in all of the determinant flow variants. 

%
%

%
%

\bibliographystyle{unsrt}
\bibliography{det}

\end{document}